
\input amstex
\documentstyle{ams-j}
\NoBlackBoxes
\nologo
\topmatter
\title
What is the Rees algebra of a module?
\endtitle
\rightheadtext{Rees algebra of a module}
\author
David Eisenbud, Craig Huneke and Bernd Ulrich
\endauthor
\address
Mathematical Sciences Research Institute,
1000 Centennial Dr.,
Berkeley, CA 94720
\endaddress
\email
de\@msri.org
\endemail
\address
Department of Mathematics, University of Kansas,
Lawrence, KS 66045
\endaddress
\email   
huneke\@math.ukans.edu
\endemail
\address
Department of Mathematics,
Michigan State University,
East Lansing, MI 48824
\endaddress

\email
ulrich\@math.msu.edu
\endemail
\thanks
All three authors were partially supported by the NSF.
This paper was last revised August 29, 2001.
\endthanks

\def\R{\Cal R}

\def\ZZ{\Bbb Z}

\def\th{{\text th}}

\def\m{\bold m}
\def\ra{\rightarrow}



\def\Hom{\hbox{\rm Hom}}

\def\Sym{\hbox{\rm Sym}}
\def\Im{\hbox{\rm Im}}

\abstract
In this paper we show that the Rees algebra can be made into a functor
on modules over a ring in a way that extends its classical definition
for ideals.  The Rees algebra of a module $M$ may be computed in terms
of a ``maximal'' map $f$ from $M$ to a free module as the image of the
map induced by $f$ on symmetric algebras.  We show that the analytic
spread and reductions of $M$ can be determined from any embedding of
$M$ into a free module, and in characteristic 0---but not in
positive characteristic!---the Rees algebra
itself can be computed from any such embedding.
\endabstract
\endtopmatter
\document

The Rees algebra of an ideal $I$ in a ring $R$, namely
$
\R=\oplus_{n=0}^\infty I^n=R[It]\subset R[t],
$
plays a major role in commutative algebra and in
algebraic geometry since 
$Proj(\R)$ is the blowup of $Spec (R)$
along the subscheme defined by $I$.
Several authors have found it useful to generalize this
construction from ideals to modules; see for instance
Gaffney and Kleiman [1999],
Katz [1995], 
Katz and Kodiyalam [1997], 
Kleiman and Thorup [1997],
Kodiyalam [1995], 
Liu [1998],
Rees [1987],
Simis, Ulrich, and Vasconcelos [1998, 2000],
and Vasconcelos [1994],
who define the Rees algebra of a module satisfying
one or another hypothesis. Usually this hypothesis was
tailored to approach the problem(s) the authors were
interested in solving.

The goal of this paper is to 
clarify the definition for arbitrary finitely generated 
modules over a Noetherian ring $R$. Our interest in this clarification arose
through our work on generalized prinicipal ideal theorems
and the heights of ideals of minors, where we heavily use
Rees algebras (see Eisenbud-Huneke-Ulrich [2000a-b, 2001]).
It seems worthwhile to understand the differences
and similarities of the various approaches from the papers above, and
to make the definition as functorial as possible.
Even for ideals there is
a problem: in the depth 0 case
it is not clear from the definition above whether
the Rees algebra depends on the embedding of $I$ in
$R$. 

A natural approach is to define the Rees algebra of a module as the
symmetric algebra of the module modulo torsion (that is, modulo
elements killed by nonzerodivisors of the ring). This does {\it not\/}
give a satisfactory definition in all cases in the sense that it may give the wrong answer
even for an ideal, if the ring is
not a domain. In general, as was well-known, it is a good definition when the module
$M$ ``has a rank'', i.e., when $M$ is free of constant rank
locally at the associated primes of $R$. This hypothesis is sufficient
for many applications; however, for example, it is not necessarily preserved replacing
$M$ by $M/xM$ and $R$ by $R/xR$ even if $x$ is a non-zerodivisor on $M$ and $R$. 

Another alternative is to consider a module $M$ together with an
embedding into a free module $G$, and define the Rees algebra of $M$
to be the subalgebra of the symmetric algebra of $G$ generated by $M$;
more generally, for any map $g: M\to N$ we define $\R(g)$, the Rees
algebra of $g$, to be the graded $R$-algebra which is the image of the
map $\Sym(g): \Sym(M)\to \Sym(N)$. One may then try to define the Rees
algebra of $M$ as $\R(g)$ for an embedding $g$ of $M$ into a free
module $G$.

However, the result may depend on the chosen embedding  $g$.
In Section 1 we give as example a principal ideal $I$ in an
Artinian ring of characteristic $p>0$, and an embedding $g: I\to R^2$
such that $\R(g)$ is not isomorphic to $\R(I)=\oplus_{n= 0}^\infty I^n$.
Of course $\R(I)$ may also be expressed as $\R(i)$, where $i$ is the
inclusion of $I$ as an ideal of $R$.  Thus $\R(g)$ depends on $g$, not
only on $M$. 

Here we take a third alternative:
\medskip
\noindent{\bf Definition 0.1.} If $R$ is a ring and $M$
is an $R$-module, we define the {\it Rees algebra\/} of $M$ to 
be 
$$
\R(M)=\Sym(M)/(\cap_gL_g)
$$
where the intersection is taken over all maps $g$ from $M$ to
free $R$-modules, and $L_g$ denotes the kernel of $\Sym(g)$.
\medskip

Although the definition may at first appear somewhat complicated, it is at least
obviously functorial: if $h:M\to N$ is a homomorphism of $R$-modules
then for every map from $N$ to a free module $g:N\to G$ the map $gh$
is a map from $M$ to a free module, so $\Sym(h):\Sym(M)\to\Sym(N)$
induces a map $\R(M)\to\R(N)$.  As the symmetric algebra functor
preserves epimorphisms, so does the Rees algebra functor. 

In Section 1 we solve the problem of computing $\R(M)$ by showing that
(for finitely generated $M$) $\R(M)=\R(f)$ for any map from $M$ to a
free module $F$ such that the dual map $F^*\to M^*$ is
surjective. This implies that forming the Rees algebra of a finitely
generated module over a Noetherian ring commutes with flat base
change.  We also show that if $g:M\to R$ is an embedding in a free
module of rank 1, then $\R(M)=\R(g)$, so that $\R(M)$ agrees with the
classical definition for ideals (in particular, this shows that the
classical definition is independent of the choice of representation of
$M$ as an ideal). Moreover, we show that in many cases $\R(M)$ can be
computed from any embedding.  The following is a special case of what
we prove:

\proclaim{Theorem 0.2} Let $R$ be a Noetherian ring
and let $M$ be a finitely generated $R$-module. 
If $R$ is torsion free over $\ZZ$,
or $R$ is unmixed and generically
Gorenstein, or $M$ is free locally at
each associated prime of $R$, 
then $\R(M)\cong\R(g)$
for any embedding $g: M\to G$ of $M$ into a free module $G$.
\endproclaim 
\goodbreak

In Section 2 we use our Rees algebra construction to introduce
analytic spread and integral dependence for arbitrary modules.
We prove that for any embedding $g$ of $M$ into
a free module, the natural map $\R(M)\to \R(g)$  has nilpotent
kernel. It follows that analytic spread and integral
dependence can be computed in $\R(g)$.

In two future papers [2000a, 2000b] we will apply the notions
developed here to obtain new generalized principal ideal theorems and
results on heights of ideals of minors of a matrix.

\goodbreak	
\bigskip
\bigskip
\centerline{\bf 1. Rees Algebras\rm}
\bigskip

We begin with an example showing that the Rees algebra of a module
cannot be defined from an arbitrary embedding in a free module, even
when the module is an ideal:
\medskip

\noindent{\bf Example 1.1.} Let $k$ be a field of characteristic $p$, let 
$$R= k[X,Y,Z]/((X^p,Y^p) + (X,Y,Z)^{p+1}),
$$
write $x,y,z$ for the images of
$X,Y,Z$ in $R$, and take $M$ to be the ideal
$M = Rz \cong R/(x,y,z)^p$.
Write $g_1:M\to R$ for the inclusion $M=Rz\subset R$, and let
$g_2: M\to R^2=Rt_1\oplus Rt_2$ be the map sending $z$ to
$xt_1+yt_2$. It is easy to see that $g_2$ is also an embedding.
The algebra $\R(g_1)$ is the same as the
classical Rees algebra $\R=\oplus_{n=0}^\infty (z^n)$, and has
$p^\th$ graded component $(z^p)\neq 0$. On the other hand
$$
\R(g_2)_p = R(x^p t^p_1 + y^p t^p_2) = 0,
$$
and it follows that $\R(g_2)$ cannot surject onto the classical
$\R=\R(g_1)$ by any graded homomorphism, so
$\R(g_2)\not\cong \R(M)$ as graded rings. (Computation using
the computer algebra system Macaulay2 shows
that (at least in low characteristics)
the vector space dimension of
$\R(M)$ is 33, while that of $\R(g_2)$ is 42,
so they are not abstractly isomorphic either.)
\medskip

To compute the Rees algebra of a module, we use the
following notion.
\medskip

\noindent{\bf Definition 1.2.} Let $R$ be a ring and
let  $M$ be an $R$-module. We say that a map $f$ from $M$
to a free $R$-module $F$ is
{\it versal} 
if every map from $M$ to a free module factors through $f$.
\medskip

It follows at once from the definition that if $f:M\to F$ is versal 
then the kernel of $\Sym(f)$ maps to zero in $\R(M)$; and
thus if $f:M\to F$ is a versal map to a free $R$-module $F$, then
$\R(M)=\R(f)$. With a finiteness assumption
it is easy to find such a map:

\proclaim{Proposition 1.3} Let $R$ be a ring, let $M$ be a finitely
generated $R$-module, and let $f:M\to F$ be a map from $M$ to 
a free $R$-module $F$.
If the dual map $F^*\to M^*$ is an epimorphism 
then $f$ is versal and $\R(M)=\R(f)$. In particular,
formation of the the Rees algebra of a finitely
generated module over a Noetherian ring commutes with
flat base change.
\endproclaim

\demo{Proof} Let $g:M\to G$ be a map from $M$ to a free
$R$-module. We must show that $g$ factors through $f$.
Since $M$ is finitely generated, $g$ factors through
a finitely generated free summand of $G$, and we may
assume that $G$ is finitely generated. If follows that
the dual $G^*$ is free. Consequently we may write $g^*=f^*h$
and since $F$ and $G$ are reflexive the desired factorization is
$g=h^{*}f$. The last statement
follows from the argument given just before the Proposition.

If $R$ is Noetherian and $M$ is a finitely generated
$R$-module, then the dual $f$ of any epimorphism
from a finitely generated free module onto $\Hom_R(M,R)=M^*$ will satisfy
the hypothesis of the Proposition. 
If now $S$ is a flat $R$-algebra and $f$ is such a map, then because
$M$ is finitely presented $S\otimes_R\Hom_R(M,R)=\Hom_S(M,S)$, so
$\R(S\otimes_RM)=\R(S\otimes_R f)=S\otimes_R\R(M)$ as required.
\qed\enddemo

\medskip

\noindent{\bf Example 1.1 continued.} 
By Proposition 1.3 the map $g_2$ is not versal. It is easy to check that
$M^*$ requires 3 generators and a
versal map from $M$ to $R^3$ may be written as 
$$
f:M=Rz \ra R^3 = Rt_1 \oplus Rt_2 \oplus Rt_3;\quad z\mapsto xt_1+yt_2+zt_3.
$$
We have
$$
\R(f)_p = R(x^p t_1^p + y^p t_2^p  + z^p t_3^p) = Rz^pt^p_3 = \R(M)_p,
$$
as implied by Proposition 1.3.
\goodbreak
\medskip

Next we show that our definition of the Rees algebra agrees
with the classical notion for ideals, which is thus
independent of the embedding of the ideal in $R$.

\proclaim{Theorem 1.4} Let $R$ be a ring, let $M$ be a finitely generated
ideal of $R$, and let $g:M\to R$ be the inclusion map.
The natural map  from $\R(M)$
to the classical Rees algebra $\R=\oplus_{n=0}^\infty I^n=\R(g)$
is an isomorphism.
\endproclaim

\demo{Proof} Let $f:M\to F$ be a versal map to a 
free module $F$, of rank $n$, say. 
Since $f$ is versal we may
find a map $h:F\to R$ so that $hf=g$. We must show
that $\phi:=\Sym(h)$ is a monomorphism on the subring
of $\Sym(F)$ generated by $f(M)$.

Write $\Sym(F)=R[t_1,\dots,t_n]$, and $\Sym(R)=R[z]$.
Let $m_1,\dots,m_s$ be generators of $M$, and
write $a_i=g(m_i)\in R$, so that $h(f(m_i))=a_iz$.
Let $S=R[x_1,\dots,x_s]$ be a polynomial ring,
and consider the map $\psi: S\to Sym(F)$ sending
$x_i$ to $f(m_i)$. We must show that the kernel of 
$\psi$ is the same as the kernel of $\phi\psi$.
Giving each $x_i$ degree 1,
the kernel of $\phi\psi$ is homogeneous, so it
suffices to show that if
$u\in S$ is a form of degree $d$ such that $\phi\psi(u)=0$
then $\phi\psi=0$.
We do induction on $d$,
the case $d=0$ being obvious.

We may write
${u = \sum^s_{i=1} u_i x_i}$ where the $u_i$ are forms of degree $d-1$.
We see that 
$$
0
=\phi\psi(u)
=\phi\psi(\sum_i(x_iu_i))
=\phi\psi(\sum_i(a_iu_i))z,
$$
so $\phi\psi(\sum_i(a_iu_i))=0$. By our induction hypothesis,
$\psi(\sum_i(a_iu_i))=0$ too.

We may expand each $\psi(u_i)$ in the form 
$\psi(u_i)=\sum_\alpha r_{i,\alpha}t^\alpha$,
where the sum runs over all multi-indices $\alpha$ of degree $d-1$,
and thus $\sum_i\sum_\alpha a_ir_{i,\alpha}t^\alpha = 0$.
Since the distinct monomials $t^\alpha$ are linearly independent,
we have $\sum_ia_ir_{i,\alpha}=0$ for each $\alpha$.
By our hypothesis that $f$ is an embedding, the $f(m_i)$
satisfy the same linear relations as the $a_i$, so we get
$\sum_if(m_i)r_{i,\alpha}=0$ for each $\alpha$, and finally
$\psi(u)=\sum_i\sum_\alpha f(m_i)r_{i,\alpha}t^\alpha = 0$,
as required.
\qed\enddemo

Note that a versal map from $M$ to a free module has the
same image as the natural map from $M$ to its double dual. This image
is called the {\it torsionless quotient\/} of $M$. Any map from
$M$ to a free module factors uniquely through the torsionless quotient
of $M$. The following result gives conditions under
which the Rees algebra of a torsionless module can be deduced from
any inclusion into a free module. For convenience in applications
we state it without the torsionless hypothesis.

\proclaim{Theorem 1.5}  Let $R$ be a Noetherian ring, let $M$ be a
finitely generated $R$-module, and let $g: M \to G$ be a map to
a free $R$-module $G$ inducing an inclusion on the torsionless quotient of $M$.
If for 
each associated prime $Q$ of $R$ either 
$R_Q$ is Gorenstein, 
or $M_Q$ is free,
or $R_Q$ is ${\Bbb Z}$-torsion free, 
then the natural epimorphism $\R(M) \to \R(g)$ is an isomorphism.
\endproclaim

\noindent{\bf Proof.} Replacing $M$ by its torsionless quotient,
we may assume that $g:M\to G$ is an inclusion.

Let $f:M \to F$ be a versal map from $M$ to a free module,
and let $\R=\R(M)=\R(f)$ be
the Rees algebra of $M$.
Let $h: F\to G$ be a map with $g=hf$ and let $\phi: \R\to\R(g)$
be the induced epimorphism. 
To prove the injectivity of $\phi$ we may
replace $R$ by $R_Q$, where $Q$ is an associated prime of $R$.  

If $R$ is Gorenstein (and hence Artinian), then free modules are
injective, and therefore any monomorphism from a module to a free
module is versal.  Similarly, if $M$ is free, then since $R$ has depth
0 any monomorphism from $M$ to a free module splits, and thus again is
versal. In either case we see that $f$ and $g$ are both versal, and
the injectivity of $\phi$ follows from the functoriality of the Rees
algebra.

Finally, we treat the case where $R$ is ${\Bbb Z}$-torsion free.
It suffices to consider the case
where $F = G \oplus H$ with $H$ free, and $h$ is the natural projection.
Let $\phi$ be the natural epimorphism $\R\to\R(g)$.
Set $J=H\cdot \Sym(F)$, the kernel of $\Sym(h)$.
Since $g$ is an inclusion,  $\phi$ is an
injection in degree 1, and we must show that
$\phi$ is an injection in every degree, or equivalently $J\cap \R=0$.
This follows from the next Lemma:

\proclaim{Lemma 1.6}  Let $R$ be a
ring  and let $F$ be a free $R$-module.
Let $M$ be a submodule of $F$, and let $\R$ be the subalgebra of
$\Sym(F)$ generated by $M$. Let $H$ be a summand of $F$ with
$H\cap M=0$. If $d!$ is a nonzerodivisor in $R$ then
$(H\cdot\Sym(F))\cap\R_d=0$.
\endproclaim

\demo{Proof} It is enough to prove that the
Lemma holds after localizing at each maximal ideal of $R$,
so we may assume that $R$ is local.
Let $t_1,\dots,t_n$ be a basis of $F$. Since $R$
is local we may suppose that $H$ is generated by
$t_{m+1},\dots,t_n$.  Set $J=H\cdot \Sym(F)$.  We may assume that $d >
1$.  Since $(d-1)!$ is a non zerodivisor as well, we know by induction
that $J \cap \R_{d-1} = 0$.  Writing $\partial_i =
\frac{\partial}{\partial t_i}$, one has $\partial_i (\R) \subset \R$
for every $i$ because the $R$-algebra $\R$ is generated by linear
forms, $\partial_i(J) \subset J$ for every $i \le m$, and
$\partial_i(J^2) \subset J$ for every $i$.

Now if $u \in J \cap \R_d$ then $\partial_i (u) \in J \cap\R_{d-1} = 0$ 
for every $i \le m$.  Since $d!$ is a non zerodivisor
on $R$ it follows that $u \in R[t_{m+1}, \dots, t_n]_d
\subset J^2$.  Thus $\partial_i(u) \in J \cap \R_{d-1} = 0$ for every
$i$, and hence $u = 0$ because $d!$ is a nonzerodivisor.
\qed\enddemo

To connect our definition of the Rees algebra with the torsion in the
symmetric algebra, suppose that $R$ be a Noetherian ring, $M$ is a
finitely generated $R$-module, and $\Cal A$ is the $R$-torsion of
Sym$(M)$.  If $M_Q$ is free for every associated prime $Q$ of $R$ and
if $f:M\to F$ is a versal map to a free module, then Sym$(f)$
induces an isomorphism Sym$(M)/{\Cal A} \tilde{\ra} \R(M)$. This is
because Sym$(f)_Q$ is injective for every associated prime $Q$ of $R$.

If we are only interested in the reduced structure of the 
Rees algebra of $M$, then we can compute it from any 
embedding of the torsionless quotient of $M$. To
prove this we need to identify the minimal primes of $\R(M)$.

\proclaim{Proposition 1.7} Let $R$ be a Noetherian ring and let $M$ be a
finitely generated $R$-module. There is a one-to-one correspondence
between the minimal primes of $\R(M)$ and the minimal primes
of $R$ given by $P\mapsto R\cap P$, and similarly for associated 
primes.
\endproclaim

\demo{Proof}
Notice that $\R(M)$ is an $R$-subalgebra of a polynomial
ring on finitely many variables over $R$. Any minimal prime of
$\R(M)$ is thus a contraction of a minimal prime of the 
polynomial ring, and the distinct minimal primes of 
the polynomial ring contract to the distinct minimal primes of $R$.
The proof for associated primes is similar.
\qed\enddemo

\proclaim{Proposition 1.8}  
Let $R$ be a Noetherian ring, let $M$ be a
finitely generated $R$-module,
and let $g: M \to G$ be a map to
a free $R$-module $G$ inducing an inclusion on the 
torsionless quotient of $M$.
The kernel of the natural epimorphism $\R(M)\to \R(g)$
is nilpotent.
\endproclaim

\demo{Proof}
By Proposition 1.7 every minimal
prime of $\R(M)$ contracts to a minimal prime of $R$,
so  it is enough to 
prove the result after localizing at a minimal
prime of $R$. Thus we may assume that $(R,\m)$ is Artinian
and local. We may also replace $M$ by its torsionless quotient
and assume that $g$ is a monomorphism.

Let $f:M\to F$ be a versal map from $M$ to a free module,
and let $K$ be the kernel of 
the natural epimorphism $\R(M)=\R(f)\to \R(g)$.
Suppose first that $M$ has no free summand. Because $R$
is local, it follows that $\Im f\subset \m F$. By the
versality of $f$, the map
$g:M\to G$ factors as $g=hf$. The kernel $K$
must be contained in the 
positive degree part of $\R(M)$, which is contained
in $\m\Sym(F)$. As $\m$ is nilpotent, $K$ is nilpotent as well.

In the general case, let $H$ be a maximal free submodule of $M$.
Since any inclusion of finitely generated
free modules over an Artinian ring splits,
we may write  $G=G'\oplus H$ and $M=M'\oplus H$
in such a way that $g=g'\oplus 1_H$. The map
$\R(M)\to \Sym(G)$ is obtained from the map
$\R(M')\to \Sym(G')$ induced by $g'$
by adjoining polynomial variables. As the kernel
of the latter map is nilpotent, the desired result
follows.
\qed\enddemo

\goodbreak
\bigskip
\bigskip
\centerline{\bf 2. Integral Dependence\rm}
\bigskip

In this section we introduce general definitions of integral
dependence and of analytic spread for modules that we will apply
elsewhere. 
\medskip

\noindent{\bf Definitions 2.1.}\quad  Let $R$ be a ring, 
let $M$ be a finitely generated
$R$-module, and $U \subset L$ submodules of $M$.
\roster
\item Let $U',L'$ be the images of $U,L$ in ${\R}(M)$ and consider
the subalgebras $R[U'] \subset R[L']\subset {\R}(M)$.  We say $L$ is {\it
integral} over $U$ {\it in} $M$ if the ring extension $R[U'] \subset R[L']$
is integral.

\item We say $M$ is {\it integral} over $U$ or $U$ is a {\it reduction} of
$M$, if $M$ is integral over $U$ in $M$.
\endroster
In the situation of (1) the Rees algebra of $L$ maps to the Rees algebra
of $M$, so if $L$ is integral over $U$, then $L$ is integral over $U$ in $M$.

\goodbreak
\proclaim{Theorem 2.2}  Let $R$ be a Noetherian ring, let $M$ be a
finitely generated $R$-module, $U \subset L$ submodules of $M$, let $f:M
\ra F$ be a versal map from $M$ to a free $R$-module.  The following are
equivalent:
\roster
\item$L$ is integral over $U$ in $M$.
\item For every minimal prime $Q$ of $R$, the module
$L'$ is integral over $U'$ in
$M'$, where $'$ denotes images in $F/QF$.
\item For every map
$M\to G$ to a free $R$-module and
for every homomorphism $R \ra S$ to a domain $S$, 
the module $L'$ is integral
over $U'$ in $M'$, where $'$ denotes tensoring with $S$ and taking images in
$S\otimes_RG$.
\item For every homomorphism $R \ra V$ to a rank one discrete valuation
ring $V$ whose kernel is a minimal prime of $R$, we have $U' = L'$, where $'$
denotes tensoring with $V$ and taking images in $V\otimes_RF$.
\item (Valuative Criterion of Integrality) 
For every map
$M\to G$ to a free $R$-module and 
every homomorphism $R \ra V$ to a rank one discrete valuation
ring $V$, we have $U' = L'$, where $'$ denotes tensoring with $V$ and taking images
in $V\otimes_RG$.
\endroster
\endproclaim
\goodbreak

\demo{Proof} 
By the functoriality of the Rees algebra, 
we may assume $G=F$ in parts (3)
and (5). As ${\R}(M)$ embeds into ${\R}(F) = \text{Sym(F)}$,
we may replace $M$ by $F$ in (1).  In Items (2)--(5), the rings $R/Q, S,
V$ are domains, and hence by Theorem 1.5 the embedding of $M'$ into
the free modules $F/QF, F\otimes_RS, F\otimes_RV$, 
respectively, can be used to define
${\R}(M')$.  Thus we may replace $M$ by $F$ in these items as well.

Now it is obvious that (1) implies (3).  Part (1) follows from (2) since
${\R}(F)/\sqrt{0} \subset \prod_Q {\R}(F/QF)$, where $Q$ ranges over
all minimal primes of $R$ and $F/QF$ is considered as a module over $R/Q$.
Finally, the equivalence of (2) and (4) and of (3) and (5) has been shown
in Rees [1987, 1.5(ii)].
\qed
\enddemo

We see from Theorem 2.2 that our definition of integrality differs
from that of Rees [1987, p. 435] when $R$ is not a domain. 
Rees' definition amounts to saying that
for every minimal prime 
$Q$ of $R$, the module $L'$ is integral (in our sense or his) 
over $U'$ in $M'$, where now
$'$ denotes images in $M/QM$. If, for example,
$k$ is a field, $R = k[x]/(x^2)$, and $U = 0 \subset L =
M = (x)/(x^2)$, then $M$ is integral over $U$ in our sense but
not in the sense of Rees.

\medskip
\noindent{\bf Definition 2.3.}
Let $R$ be a Noetherian local ring with residue
field $k$ and let $M$ be a finitely generated $R$-module.  The {\it
analytic spread} $\ell(M)$ of $M$ is the Krull dimension of $k \otimes_R
{\R}(M)$.

\medskip
By way of illustration, we remark that the analytic spread
of a finitely generated module $M$ 
over an Artinian local ring $(R,\m, k)$ 
is equal to the rank $r$ of
a maximal free summand $H$ of $M$ (and is also equal
to the dimension of the Rees algebra of the module). This is
because any homomorphism of $M/H$ to a free module $F$
carries $M/H$ into $\m F$, which generates a nilpotent ideal
of $\Sym(F)$. Thus $(\R(M)/\m \R(M))_{\text{red}}=\R(M)_{\text {red}}$
is a polynomial ring over $k$ on $r$ variables. 

If $k$ is infinite one can show as in the case of ideals, using a
homogeneous Noether normalization of $k \otimes_R {\R}(M)$ and
Nakayama's Lemma, that 
$$
\ell(M) = \min\{\mu(U) \mid U \hbox{ is a reduction of }M\}.
$$
Furthermore, $\ell(M) \leq \mu(M)$ and equality holds if and only
if $M$ has no proper reduction.

\proclaim{Proposition 2.4}
Let $R$ be a Noetherian local ring with residue field $k$,
let $M$ be a finitely generated $R$-module,
and let $g: M\to G$ be a map to a free $R$-module $G$.
If $g$ induces an inclusion on the torsionless quotient of $M$, 
then
$\ell(M)=\dim k\otimes_R\R(g)$.
\endproclaim

\demo{Proof} Proposition 1.8 shows that $\R(M)$ differs from 
$\R(g)$ only by a nilpotent ideal, and thus the same
holds after tensoring with $k$.
\qed\enddemo
\goodbreak

\centerline{\bf Bibliography}
\bigskip
\refstyle{A}

\enddocument